\documentclass[12pt]{article}

\usepackage{amsmath, amsbsy, amsthm, amsfonts, array,graphics}

 \DeclareMathOperator{\der}{d}
 \DeclareMathOperator{\Vol}{Vol}
 
 \DeclareMathOperator{\B}{B}
 
 \DeclareMathOperator{\Ric}{Ric}
 \DeclareMathOperator{\Rm}{Rm}

 \newtheorem{thm}{Theorem}[section]
 \newtheorem{cor}[thm]{Corollary}
 \newtheorem{lem}[thm]{Lemma}

 \numberwithin{equation}{section}

\title
{\textbf{ \normalsize{Local Volume Estimate for Manifolds with $L^2$-bounded Curvature}}}
\author
{\normalsize{Ye Li}\\
\it{\footnotesize{Department of Mathematics,}} \\
\it{\footnotesize{Princeton University, 08544, U.S.A.}}\\
\it{\footnotesize{email address: yeli@math.princeton.edu}} \\
}
\date{}

\begin{document}
\maketitle
\begin{center}
\begin{minipage}{12cm}
    \hspace{0.5cm}\textbf{Abstract:}
 we obtain a local volume growth for complete, noncompact Riemannian manifolds with small integral bounds and with
Bach tensor having finite $L^2$ norm in dimension 4.
\end{minipage}
\end{center}

\vskip 0.2cm
\begin{center}
\begin{minipage}{12cm}

\end{minipage}

\end{center}
\vskip 0.2cm

\section{Introduction}

\ \ \ It is important to study asymptotic behavior of complete
manifold without the assumption of pointwise Ricci curvature
bound. A volume growth
and curvature decay result was obtained in \cite{3} for various classes of complete,
noncompact, Bach-flat metrics in dimension 4.  Some similar results were also claimed in \cite{1}.

In this note we consider a more general case, that is, the Bach tensor may not necessarily vanish. Since Bach tensor can be viewed as
a second derivative of the Ricci tensor,  there will be a priori no $L^{p}$ bound for it, where $p > 2$. So we may consider the
case that the $L^2$ norm of the Bach tensor is finite.  Our main result is to give a local volume estimate:

\begin{thm}\label{thm1.1}
Let $X$ be a complete, noncompact 4-dimensional Riemannian manifold. Let $B(p,r)$ be a geodesic ball around the point $p$.
Assume that there holds the following local Sobolev inequality: for any open subset $\Omega$,
 $$\|f\|^2_{L^4(\Omega)}\le C_s(\Omega)\|\nabla f\|^2_2,\ \forall f\in C_0^{\infty}(\Omega).$$
Then there exist constants $\varepsilon_0$ and $C$ (depending on the Sobolev constant $C_s(B(p,r))$) such that if
$$\|\Rm\|_{L^2(B(p,2r))}\le \varepsilon_0$$
and $\B\in L^2(B(p,2r))$, then
$$\Vol(B(p,r))\le C r^4.$$
\end{thm}

If the Bach tensor does not vanish, then a direct computation shows that
\begin{equation*}
\triangle \Ric=\Rm*\Ric+\B,
\end{equation*}
where $\B$ is the Bach tensor.  A standard argument to obtain the
bound for the Ricci tensor is to use the elliptic Moser iteration
for this equation. However, as we mentioned above, we can't assume
$\B\in L^{p}$ for $p>2$, for this will automatically give the
regularity for the Ricci tensor. So it is not obvious to apply the
elliptic Moser iteration directly, since now we consider an
inhomogeneous equation.

To overcome this difficulty, we use the Ricci flow to smooth the Riemannian metric, which was first considered by Bemelmans, Min-Oo and Ruh \cite{2}.
Notice that since we only consider the local case, what matters is not the global $L^2$ bound on curvature but the local bound, that is, the $L^2$
norm of curvature on each geodesic ball of fixed radius. Also, a global heat flow will not control such a local bound.
So, instead, we will use local Ricci flow , which was first used by D. Yang \cite{4}.  In that paper, a simple form of Moser iteration was
applied to a local nonlinear heat equation.  And we found that this argument works in our settings to obtain a pointwise local bound for the
curvature tensor of the regularized metric via the local Ricci flow.

We end the introduction with a brief outline of the note.  In Section 2, we will prove the Moser iteration for the local heat flow.
The local existence of the Ricci flow will be discussed in Section 3. And the local bound for the curvature tensor of the regularized
 metric will be obtained in Section 4.  Finally Theorem \ref{thm1.1} will be proved in Section 5.

\textbf{Acknowledgement.}  The author would like to thank his advisor, Professor Gang Tian, for many helpful and
stimulating discussions and for bringing his attention to the paper \cite{3}.

\vskip 1cm

\section{Moser Iteration for a Local Heat Flow}

\ \ \ Fix an open set $B_0\subset X$ and a smooth compactly supported function $\phi\in C_0^{\infty}(B_0)$.

Let $g(t),0\le t \le T$, be a 1-parameter family of smooth Riemannian metrics.
Let $\nabla$ denote the covariant differentiation with respect to the metric $g(t)$ and
$-\triangle$ be the corresponding Laplace-Beltrami operator.
Let $A>0$ be a constant that satisfies the standard Sobolev inequality
$$\left(\int_{B_0}f^{4}\der V_g\right)^{\frac 12}\le A \int_{B_0}|\nabla f|^2\der V_g,\ f\in C_0^{\infty}(B_0),$$
with respect to each metric $g(t), 0\le t\le T$.

Assume that for each $t\in [0,T]$,
$$\frac{1}{2}g_{ij}(0)\le g_{ij}(t)\le 2 g_{ij}(0)\ \ \text{on}\ B_0.$$
All geodesic balls in this section are defined with respect to the metric $g(0)$,
and therefore, are fixed open subsets of $X$, independent of $t$.

We want to study the following heat equation:
\begin{align}\label{1}
\frac{\partial f}{\partial t}\le \phi^2(\triangle f+uf)+2a\phi|\nabla \phi||\nabla f|+b(|\nabla \phi|^2-\phi\triangle \phi)f,\ 0\le t \le T,
\end{align}
where $f$ and $u$ are nonnegative functions on $B_0\times [0,T]$, such that
\begin{align}\label{2}
\frac{\partial}{\partial t}\der V_g\le c\phi^2u\der V_g
\end{align}
and
\begin{align}\label{3}
\left(\int_{B_0}\phi^2u^3\right)^{\frac 13}\le \mu t^{-\frac 13}.
\end{align}

The following results in this section are due to D. Yang \cite{4}.  For convenience, we give the proofs below.  Notice that our manifold is 4-dimensional.

\vskip 0.2cm
\begin{lem}\label{lem1}
Given $p>1,\  \psi\in C_0^{\infty}(B_0),\ f\in C^{\infty}(M),\ f\ge 0$,
$$\int_{B_0}|\nabla(\psi f^{\frac{p}{2}})|^2\le \frac{p^2}{2(p-1)}\int_{B_0}\psi^2f^{p-1}(-\triangle f)\der V_g
+\left(1+\frac{1}{(p-1)^2}\right)\int_{B_0}|\nabla \psi|^2f^p\der V_g.$$
\end{lem}

\noindent\textit{Proof:}  Using integration by parts, we have
\begin{align*}
\int |\nabla (\psi f^{\frac{p}{2}})|^2 =& -\int \psi f^{\frac p2} \triangle (\psi f^{\frac p2})\\
    =&\frac p2 \int \psi^2 f^{p-1} (-\triangle f)+ \int f^p|\nabla \psi|^2-\frac{p(p-2)}{4}\int \psi^2 f^{p-2}|\nabla f|^2\\
    =&\frac {p^2}{2(p-1)} \int \psi^2 f^{p-1}(-\triangle f)+\frac{p}{2(p-1)}\int \psi^2 f^{p-1} \triangle f \\
     & +\int f^p |\nabla \psi|^2-\frac{p(p-2)}{4}\int \psi^2 f^{p-2} |\nabla f|^2.
\end{align*}
On the other hand, by Cauchy inequality,
\begin{align*}
\frac{p}{2(p-1)}\int \psi^2 f^{p-1}\triangle f &= -\frac{p}{2(p-1)}\int \nabla (\psi^2 f^{p-1})\nabla f\\
    &=-\frac{p}{p-1}\int \psi f^{p-1} \nabla \psi \nabla f-\frac{p}{2}\int \psi^2 f^{p-2}|\nabla f|^2\\
    &\le \frac{1}{(p-1)^2}\int f^p |\nabla \psi|^2+\frac{p^2}{4}\int \psi^2f^{p-2}|\nabla f|^2-\frac p2 \int \psi^2 f^{p-2}|\nabla f|^2\\
    &=\frac{1}{(p-1)^2}\int f^p|\nabla \psi|^2+\frac{p(p-2)}{4}\int \psi^2 f^{p-2}|\nabla f|^2.
\end{align*}
This proves the lemma. $\Box$

\vskip 0.2cm
\begin{lem}\label{lem2}
Suppose that $f$ and $u$ are nonnegative functions on $B_0\times[0,T]$ which satisfy (\ref{1}), (\ref{2}) and (\ref{3}).
For $p\ge p' \ge p_0 >1,$ we have
\begin{align}\label{4}
\frac{\partial}{\partial t}\int \phi^{2p'}f^p+\int |\nabla(\phi^{p'+1}f^{\frac{p}{2}})|^2\le
[(p'+1)^2C\|\nabla \phi\|^2_{\infty}+C(p\mu)^3A^2t^{-1}]\int \phi^{2p'}f^p.
\end{align}
\end{lem}
\noindent\textit{Proof:} Given $p\ge p'\ge p_0>1$, we combine Lemma \ref{lem1} with (\ref{1}) and (\ref{2}) to obtain
\begin{align}\label{5}
\frac{\partial}{\partial t}\int \phi^{2p'}f^p&+2\left(1-\frac{1}{p}\right)^2\int|\nabla(\phi^{p'+1}f^{\frac p2})|^2 \nonumber \\
    \le &\int [\phi^2(\triangle (\phi^{2p'}f^p)+C\phi^{2p'}uf^p)+2a\phi|\nabla \phi||\nabla (\phi^{2p'}f^p)| \nonumber \\
      &+b(|\nabla \phi|^2-\phi\triangle \phi)\phi^{2p'}f^p].
\end{align}
Now we estimate each term of the right hand side.
\begin{align*}
\int \phi^2 \triangle (\phi^{2p'}f^p)&=-2\int \phi \nabla \phi \nabla (\phi^{2p'}f^p)\\
   &=-\int 2\phi \nabla \phi (2p'\phi^{2p'-1}\nabla \phi \cdotp f^p+p\phi^{2p'}f^{p-1}\nabla f)\\
   &=-4p'\int \phi^{2p'}f^p|\nabla \phi|^2 -2p\int \phi^{2p'+1}f^{p-1}\nabla\phi\nabla f\\
   &\le 4p'\int \phi^{2p'}f^p|\nabla \phi|^2+2p\left(\int \phi^{2p'}|\nabla \phi|^2 f^p\right)^{\frac 12}
               \left(\int\phi^{2p'+2}f^{p-2}|\nabla f|^2\right)^{\frac{1}{2}}\\
   &\le 4p'\int \phi^{2p'}f^p|\nabla \phi|^2+\frac{pC}{\varepsilon} \int \phi^{2p'}|\nabla \phi|^2 f^p
         +p\ \varepsilon \int\phi^{2p'+2}f^{p-2}|\nabla f|^2.
\end{align*}
By a similar argument the remaining terms can be estimated as follows.
\begin{align*}
\int \phi |\nabla \phi||\nabla(\phi^{2p'}f^p)|\le & 2 p'\int\phi^{2p'}f^p|\nabla \phi|^2+p\int\phi^{2p'+1}f^{p-1}|\nabla \phi||\nabla f|,\\
-\int (\phi\triangle \phi)\phi^{2p'}f^p=&\int \nabla \phi \nabla (\phi^{2p'+1}f^p) \\
            = & (2p'+1)\int \phi^{2p'}|\nabla \phi|^2f^p+p\int\phi^{2p'+1}f^{p-1}\nabla \phi \nabla f,\\
-\int \phi^{2p'+2}f^{p-1}\triangle f=& \int \nabla (\phi^{2p'+2}f^{p-1})\nabla f\\
            =& (2p'+1)\int \phi^{2p'+1}f^{p-1}\nabla \phi \nabla f+(p-1)\int \phi^{2p'+2}f^{p-2}|\nabla f|^2.
\end{align*}
So it is easy to see that each term of right hand side of (\ref{5}) has the form of
$$\int |\nabla \phi|^2\phi^{2p'}f^p,\ \int \phi^{2p'+2}f^{p-2}|\nabla f|^2\ \text{or}\ \int u\phi^{2p'+2}f^p,$$
where
\begin{align*}
\int \phi^{2p'+2}f^{p-2}|\nabla f|^2=&\left(\frac 2p\right)^2\int |\phi^{p'+1}\nabla f^{\frac p2}|^2\\
        \le & 2 \left(\frac 2p\right)^2 \varepsilon \int |\nabla(\phi^{p'+1}f^{\frac p2})|^2+
                2 \left(\frac 2p\right)^2 (p'+1)^2 C_{\varepsilon}\int \phi^{2p'}f^p|\nabla \phi|^2.
\end{align*}
Notice that if $\varepsilon$ is sufficiently small, then the first term of the right hand side can be absorbed into the
left hand side of (\ref{5}).  Therefore we have
\begin{align*}
\frac{\partial}{\partial t}\int \phi^{2p'}f^p+&2\left(1-\frac 1p\right)^2\int |\nabla(\phi^{p'+1}f^{\frac p2})|^2\\
    \le &(p'+1)^2C\int |\nabla \phi|^2\phi^{2p'}f^p+p\int u\phi^{2(p'+1)}f^p.
\end{align*}

Using H$\ddot{o}$lder, Sobolev, Cauchy inequalities, and (\ref{3}), we see that
\begin{align*}
\int u\phi^{2(p'+1)}f^p\le &\left(\int \phi^2 u^3\right)^{\frac 13}\left(\int(\phi^{2p'}f^p)\right)^{\frac 13}
                \left(\int \phi^{4(p'+1)}f^{2p}\right)^{\frac 13}\\
    \le & \mu t^{-\frac 13}\left(\int(\phi^{2p'}f^p)\right)^{\frac 13}\cdot A^{\frac 23} \left(|\nabla(\phi^{p'+1})f^{\frac p2}|^2\right)^{\frac 23}\\
    \le & (\mu t^{-\frac 13})^3\varepsilon^{-\frac 13}\int \phi^{2p'}f^p+\varepsilon^{\frac 23}A^2\int |\nabla(\phi^{p'+1}f^{\frac p2})|^2.
\end{align*}
Thus,
\begin{align*}
\frac{\partial}{\partial t}\int \phi^{2p'}f^p+&2\left(1-\frac 1p\right)^2\int|\nabla(\phi^{p'+1}f^{\frac p2})|^2\\
        \le & (p'+1)C\int |\nabla \phi|^2 \phi^{2p'}f^p+\varepsilon^{\frac 13}A^2\int |\nabla(\phi^{p'+1}f^{\frac p2})|^2
            +\varepsilon^{-\frac 13}p^3\mu^3 t^{-1}\int\phi^{2p'}f^p.
\end{align*}
Choosing $\varepsilon$ so that $\varepsilon^{\frac 23}A^2$ is sufficient small, we have
\begin{align*}
\frac{\partial}{\partial t}\int \phi^{2p'}f^p+&\int|\nabla(\phi^{p'+1}f^{\frac p2})|^2\\
    \le & [(p'+1)^2C\|\nabla \phi\|^2_{\infty}+C(p\mu)^3 A^2 t^{-1}]\int \phi^{2p'}f^p.
\end{align*}
This proves lemma \ref{lem2}.

Now given $0<\tau<\tau '<T$, let
\begin{align*}
    \psi(t)=\left\{
    \begin{array}{cc}
       0, &0\le t\le \tau,\\
       \dfrac{t-\tau}{\tau '-\tau}, &\tau\le t\le \tau ',\\
       1, &\tau '\le t\le T
    \end{array}
    \right.
\end{align*}
Multiplying (\ref{4}) by $\psi$, and noticing that $p'+1\le p^2$,
\begin{align*}
\frac{\partial}{\partial t}\left(\psi \int \phi^{2p'}f^p\right)+&\psi\int |\nabla (\phi^{p'+1}f^{\frac p2})|^2\\
    \le& [\ p^6\hat{C}(t)\psi+|\psi '| \ ]\int \phi^{2p'}f^p,
\end{align*}
where $\hat{C}(t)=C\|\nabla \phi\|^2_{\infty}+C\mu^3 A^2 t^{-1}$.
Integrating this with respect to $t$, we obtain

\vskip 0.2cm
\begin{lem}\label{lem3}
\begin{align*}
\int \phi^{2p'}f^p+\int_{\tau '}^{t}\int |\nabla (\phi^{p'+1}f^{\frac p2})|^2
    \le \left(p^6\hat{C}(\tau ')+\frac {1}{\tau '-\tau}\right)\int_{\tau}^T\int \phi^{2p'}f^p,\ \tau '\le t\le T.
\end{align*}
\end{lem}
Given $p\ge p' \ge p_0>1,\ 0\le \tau <T$, denote
$$H(p,p',\tau)=\int_{\tau}^{T}\int_{B_0}\phi^{2p'}f^p.$$

\vskip 0.2cm
\begin{lem}\label{lem4}
Given $p\ge p_0,\ 0\le \tau < \tau ' < T$,
\begin{align*}
H\left(\frac 32 p, \frac 32 p'+1, \tau '\right)\le AC [\ (\tau '-\tau)^{-1}+p^6\hat{C}(\tau ')\ ]^{\frac 32}H(p,p',\tau)^{\frac 32}.
\end{align*}
\end{lem}
\noindent\textit{Proof:} By H$\ddot{o}$lder, Sobolev inequalities,
\begin{align*}
H\left(\frac 32 p,\frac 32 p'+1, \tau '\right)=&\int_{\tau '}^{T}\int \phi^2(\phi^{2p'}f^p)^{\frac 32}\\
        \le& \int_{\tau '}^{T} \left(\int \phi^{2p'}f^p\right)^{\frac 12}\left(\int \phi^{4p'+4}f^{2p}\right)^{\frac 12}\der t\\
        \le& \left(\sup_{\tau '\le t\le T}\int \phi^{2p'}f^p\right)^{\frac 12}\ A\ \int_{\tau '}^T\int |\nabla (\phi^{p'+1}f^{\frac p2})|^2 \der t.
\end{align*}
Applying Lemma \ref{3}, we obtain the desired estimate. $\Box$

\vskip 0.2cm
\begin{thm}\label{thm1}
Let $f$ and $u$ be non-negative functions on $B_0\times[0,T]$, such that $\dfrac{\partial}{\partial t}\der V_g\le c\phi^2u\der V_g$
for some constant $c$, and
$$\frac{\partial f}{\partial t}\le \phi^2(\triangle f+uf)+2a\phi|\nabla \phi||\nabla f|+b(|\nabla \phi|^2-\phi\nabla\phi)f,\ 0\le t\le T.$$
Assume that
$$\left(\int_{B_0}\phi^2u^3\right)^{\frac 13}\le \mu t^{-\frac 13}.$$
Then given $(x,t)\in B_0\times[0,T],\ p_0>2$,
$$|\phi(x)^2f(x,t)|\le CA^{\frac {2}{p_0}}[\ \|\nabla \phi\|^2_{\infty}+t^{-1}(1+A^2\mu^3)\ ]^{\frac{3}{p_0}}
        \left(\int_0^t\int_{B_0}\phi^{2p_0-4}f^{p_0}\right)^{\frac {1}{p_0}},$$
where $C$ depends on $p_0,\ a$ and $b$.
\end{thm}
\noindent \textit{Proof: } Denote $\nu=\frac 32 and \ \eta=\nu^6$.
Fix $0<t<T$, and set
\begin{align*}
p_k^{\prime}=&\left(p_0-2\right)\nu^k+\sum_{j=0}^{k-1}\nu^j,\\
p_k=&\ p_0\nu^k,\\
\tau_k=&\ t(1-\eta^{-k}),\\
\Phi_k=&\ H(p_k,p_k^{\prime},\tau_k)^{\frac {1}{p_k}}.
\end{align*}
Applying Lemma \ref{4},
\begin{align*}
H(p_{k+1},p^{\prime}_{k+1},\tau_{k+1})\le AC \left[\ \|\nabla \phi\|^2_{\infty}+(1+\mu^3A^2)\frac{\eta}{\eta-1}t^{-1}\ \right]^{\nu}
    \eta^{k\nu}H(p_k,p_k^{\prime},\tau_{k})^{\nu}.
\end{align*}
Therefore,
$$\Phi_{k+1}\le (AC)^{\frac{\sigma_{k+1}-1}{p_0}}\left(\|\nabla \phi\|^2_{\infty}+(1+\mu^3A^2)^{\frac{\eta}{\eta-1}}t^{-1}\right)^{\frac{\sigma_k}{p_0}}
    \cdot \eta^{\frac{\sigma_{k}^{\prime}}{p_0}}H(p_0,p_0-2,0)^{\frac {1}{p_0}},$$
where $\sigma_k=\sum_{i=0}^{k}\nu^{-i},\ \sigma^{\prime}_k=\sum_{i=0}^{k}i\nu^{-i}.$
Letting $k\rightarrow \infty$, we obtain
$$|\phi^2f(x,t)|\le CA^{\frac {2}{p_0}}\left[\ \|\nabla\phi\|^2_{\infty}+t^{-1}(1+\mu^3A^2)\ \right]^{\frac{3}{p_0}}
    \left(\int_0^T\int\phi^{2p_0-4}f^{p_0}\right)^{\frac{1}{p_0}}.$$
Now let $T\rightarrow t$.  This proves the theorem.  $\Box$

\vskip 0.2cm
\begin{thm}\label{thm2}
Let $f\ge 0$.  Solve
\begin{align}\label{6}
\frac{\partial f}{\partial t}\le \phi^2(\triangle f+C_0 f^2)+2a\phi|\nabla \phi||\nabla f|+b(|\nabla \phi|^2-2\phi\triangle \phi)f,\ 0\le t\le T,
\end{align}
on $B_0\times[0,T]$.  Assume that
$$\frac{\partial}{\partial t}\der V_g\le C\phi^2f\der V_g$$
and that
$$\left(\int_{B_0}f^2_0\right)^{\frac 12}\le (5 e C_0 A)^{-1},$$
where $f_0(x)=f(x,0)$.  Then
$$|\phi^2(x)f(x,t)|\le CA(t\|\nabla \phi\|^2_{\infty}+1)^2t^{-1},$$
where $0<t<\min (T, \|\nabla \phi\|^{-2}_{\infty}),\ C=C(C_0,a,b)$.
\end{thm}
\noindent \textit{Proof: } Let $[0,T']\subset [0,T]$ be the maximal interval such that
$$e_0=\sup_{0\le t\le T'}\left(\int_{B_0}f^2\right)^{\frac 12}\le (4C_0A)^{-1}.$$
Applying Lemma \ref{1} to (\ref{6}), we have, for $0\le t\le T'$,
\begin{align*}
\frac{\partial}{\partial t}\int f^p+&2\left(1-\frac 1p\right)^2\int |\nabla(\phi f^{\frac p2})|^2 \\
    \le& p\int|\nabla \phi|^2f^p+pC_0A\left(\int f^2\right)^{\frac 12}\int |\nabla(\phi f^{\frac p2})|^2.
\end{align*}
Therefore, for $p=2$, the bound on the $L^2$ norm of $f$ implies that for $0\le t\le T'$,
$$\frac{\partial}{\partial t}\int f^2\le 2\|\nabla \phi\|_{\infty}\int f^2,$$
which implies that
$$\int f^2\le e^{2\|\nabla \phi\|^2 t}\int f_0^2.$$
In particular, if $T'< \|\nabla \phi\|^{-2}$, then
$$e_0\le e\int f_0^{\frac n2}\le (5C_0A)^{-1}<(4C_0A)^{-1}.$$
This contradicts the assumed maximality of $[0,T']$.  We can therefore assume that $T'\ge \min ((\log 2)\|\nabla \phi\|^{-2},T)$.

By the same argument of Lemma \ref{lem2}, we have an estimate of the form
$$\int f^p \cdot \int_0^t\int |\nabla (\phi f^{\frac p2})|^2\le C(t^{-1}+\|\nabla \phi\|_{\infty})\int_0^t\int f^p.$$

Therefore,
\begin{align*}
\int \phi^2 f^3\le& C (t^{-1}+\|\nabla \phi\|_{\infty})\int_0^t\int \phi^2f^3\\
    \le& C (t^{-1}+\|\nabla \phi\|_{\infty})\int_0^t\left(\int f^2\right)^{\frac 12}\left(\int (\phi f)^4\right)^{\frac 12}\der t\\
    \le& C e_0 A (t^{-1}+\|\nabla \phi\|_{\infty})\int_0^{t}\int|\nabla(\phi f)|^2\\
    \le& C e_0 A (t^{-1}+\|\nabla \phi\|_{\infty})^2 \int_0^t\int f^2\\
    \le& C At(t^{-1}+\|\nabla \phi\|_{\infty})^2e_0^3.
\end{align*}
Set
$$\mu^3=CA(1+t\|\nabla \phi\|^2_{\infty})^2 e^3$$
and notice that Theorem \ref{thm1} still holds, when $p_0\rightarrow 2$.
We then obtain the desired estimate.  $\Box$

The argument also implies the following

\begin{cor}\label{cor1}
Let $f$ satisfy the assumptions of Theorem \ref{thm2}.  Then given $u\ge 0$ such that
$$\frac{\partial u}{\partial t}\le \phi^2(\triangle u+c_0fu)+a\cdot \nabla u +bu,$$
the following estimate holds for $0\le t< \min (T, (\log 2)\|\nabla \phi\|^{-2}_{\infty})$,
$$|\phi(x)^2u(x,t)|\le CA^{\frac{2}{p_0}}(1+t\|\nabla \phi\|^2_{\infty})^2t^{-\frac{2}{p_0}}\left(\int_{B_0}u^{p_0}\right)^{\frac{1}{p_0}},$$
where $u_0(x,t)=u(x,0)$, and $C$ depends on $p_0,\ a$ and $b$.
\end{cor}

\vskip 1cm

\section{Existence of Local Ricci Flow}

\ \ \ Let $X$ be a smooth 4-manifold without boundary.
Given a smooth Riemannian metric $g_0$ and a smooth compactly supported function $\phi$,
we want to study the following evolution equation
\begin{align}\label{3.1}
    \left\{
    \begin{array}{ll}
       \dfrac{\partial g}{\partial t}&=-2\phi^2 \Ric(g),\\
       g(0)&=g_0.\\
    \end{array}
    \right.
\end{align}

\begin{thm}\label{thm3}
There exists $T>0$ such that (\ref{3.1}) has a smooth solution for $0\le t\le T$.
\end{thm}
\noindent \textit{Proof: }
Given $\varepsilon>0$, consider
\begin{align}\label{3.2}
    \left\{
    \begin{array}{ll}
       \dfrac{\partial g}{\partial t}&=-2(\varepsilon^2+\phi^2) \Ric(g),\\
       g(0)&=g_0.\\
    \end{array}
    \right.
\end{align}
We want to use DeTurck's trick so that this system can be reduced to a nonlinear,
strictly parabolic system.  First, we fix a metric $\hat{g}$ on $X$.
Let $\Gamma_{ij}^k$ and $\hat{\Gamma}_{ij}^k$ denote the Christoffel symbols of $g$ and $\hat{g}$ respectively.
Our aim is to give an expression of $\Ric(g)-\Ric(\hat{g})$.
By direct calculation, we have
$$\Gamma_{jk}^i-\hat{\Gamma}_{jk}^i=\frac{1}{2}g^{il}(g_{jl,k}+g_{kl,j}-g_{jk,l}),$$
where $g_{jl,k}=\dfrac{\partial g_{jl}}{\partial x^k}-g_{sl}\hat{\Gamma}^s_{kj}-g_{js}\hat{\Gamma}_{kl}^s$,
the covariant derivative with respect to the metric $\hat{g}$.

Recall that in local coordinates
$$R_{ijl}^p=\frac{\partial}{\partial x^i}\Gamma_{jl}^p-\frac{\partial}{\partial x^j}\Gamma^p_{il}+
    \Gamma_{ig}^p\Gamma_{jl}^q-\Gamma_{jq}^p\Gamma^q_{il},$$
and
$$R_{ik}=g^{jl}g_{hk}R^h_{jil}.$$
So
\begin{align*}
R_{ij}-\hat{R}_{ij}=&g^{kl}g_{hj}\left(\frac{\partial}{\partial x^i}\Gamma_{kl}^h-\frac{\partial}{\partial x^k}\Gamma_{il}^h\right)+\hbox{other terms}\\
        =&g^{kl}g_{hj}\left(\frac{\partial}{\partial x^i}(\Gamma_{kl}^h-\hat{\Gamma}_{kl}^h)
            -\frac{\partial}{\partial x^k}(\Gamma_{il}^h-\hat{\Gamma}_{il}^h)\right)+\hbox{other terms}\\
        =&-\frac{1}{2}g^{kl}g_{ij,kl}+\frac{1}{2}g^{kl}(g_{il,jk}+g_{jl,ik}-g_{kl,ij})+\hbox{other terms}.
\end{align*}
Now we set
$$X^p=-g^{pi}g^{kl}(g_{ik,l}-\frac{1}{2}g_{kl,i})$$
and $X=X^p\dfrac{\partial}{\partial x^p}$, then
\begin{align*}
(L_{_X} g)_{kl}=&(L_{_X} g)\left(\frac{\partial}{\partial x^k},\frac{\partial}{\partial x^l}\right)\\
        =&L_{_X}\left(g\left(\frac{\partial}{\partial x^k},\frac{\partial}{\partial x^l}\right)\right)
            -g\left(L_{_X}\frac{\partial}{\partial x^k},\frac{\partial}{\partial x^l}\right)
            -g\left(\frac{\partial}{\partial x^k},L_{_X}\frac{\partial}{\partial x^l}\right)\\
        =&X\left(g\left(\frac{\partial}{\partial x^k},\frac{\partial}{\partial x^l}\right)\right)-
                g\left(\nabla_{_X}\frac{\partial}{\partial x^k}-\nabla_{\frac{\partial}{\partial x^k}}X,\frac{\partial}{\partial x^l}\right)\\
                &-g\left(\frac{\partial}{\partial x^k},\nabla_{_X}\frac{\partial}{\partial x^l}-\nabla_{\frac{\partial}{\partial x^l}}X\right)\\
        =&g\left(\nabla_{\frac{\partial}{\partial x^k}}X,\frac{\partial}{\partial x^l} \right)+
            g\left(\frac{\partial}{\partial x^k},\nabla_{\frac{\partial}{\partial x^l}}X\right)\\
        =&\frac{\partial X^p}{\partial x^k}g\left(\frac{\partial}{\partial x^p},\frac{\partial}{\partial x^l}\right)
            +X^p g\left(\nabla_{\frac{\partial}{\partial x^k}}\frac{\partial}{\partial x^p},\frac{\partial}{\partial x^l}\right)
            +\frac{\partial X^p}{\partial x^l}g\left(\frac{\partial}{\partial x^k},\frac{\partial}{\partial x^p}\right)\\
            &+X^p\left(\frac{\partial}{\partial x^k},\nabla_{\frac{\partial}{\partial x^l}}\frac{\partial}{\partial x^p}\right)\\
        =&\frac{\partial X^p}{\partial x^k}g_{pl}+\frac{\partial X^p}{\partial x^l}g_{kp}+X^p\frac{\partial}{\partial x^p}g_{kl}.
\end{align*}
Thus
\begin{align}\label{3.3}
R_{ij}-\hat{R}_{ij}=-\frac{1}{2}g^{kl}g_{ij,kl}-\frac{1}{2}(L_{_X}g)_{ij}+\hbox{other terms}.
\end{align}
We set
$$(F(g))_{ij}=g^{kl}g_{ij,kl}+Q_{ij},$$
where $Q_{ij}$ involves the other terms in (\ref{3.3}), then
$$\Ric (g)-\Ric (\hat{g})=-\frac 12 F(g)-\frac 12 L_{_X}g.$$
We define a one-parameter diffeomorphism group $\Phi_t: X\rightarrow X$ as follows.
\begin{align*}
    \left\{
    \begin{array}{ll}
       \dfrac{\der \Phi_t(x)}{\der t}=[(\phi\circ \Phi_t^{-1})^2+\varepsilon^2]X(t,\Phi_t(x)),\\
       \ \ \Phi_0(x)=x,\\
    \end{array}
    \right.
\end{align*}
where $X=X^p\dfrac{\partial}{\partial x^p}$ given as above.

Consider the following initial value problem
\begin{align*}
    \left\{
    \begin{array}{ll}
     \dfrac{\partial \bar{g}}{\partial t}&=[(\phi\circ \Phi_t^{-1})^2+\varepsilon^2][F(\bar{g})-2\Ric(\hat{g})]-P,\\
     \bar{g}(0)&=g_0,\\
    \end{array}
    \right.
\end{align*}
where $P_{ij}=\frac{\partial}{\partial x^i}[(\phi\circ \Phi^{-1}_t)^2]\bar{g}(X,\frac{\partial}{\partial x^j})
+\frac{\partial}{\partial x^j}[(\phi\circ \Phi^{-1}_t)^2]\bar{g}(X,\frac{\partial}{\partial x^i}).$
Then a direct calculation shows that $g=\Phi_t^*(\bar{g})$ is the solution of (\ref{3.2}).  Indeed,
\begin{align*}
\frac{\partial g}{\partial t}=&\frac{\partial}{\partial t}\Phi_t^*(\bar{g})\\
    =&\Phi_t^*\left(\frac{\partial \bar{g}}{\partial t}\right)+\Phi_t^*(L_{\Phi_t}\cdot \bar{g})\\
    =&\Phi_t^*\left(\frac{\partial \bar{g}}{\partial t}+L_{[(\phi\circ \Phi_t^{-1})^2+\varepsilon^2]X}\bar{g}\right)\\
    =&\Phi_t^*\{[(\phi\circ \Phi_t^{-1})^2+\varepsilon^2][F(\bar{g})-2\Ric(\hat{g})+L_{_X}\bar{g}]\}\\
    =&-2\Phi_t^*\{[(\phi\circ \Phi_t^{-1})^2+\varepsilon^2]\Ric{\bar{g}}\}\\
    =&-2(\phi^2+\varepsilon^2)\Ric(g),
\end{align*}
where we used the following fact,
$$(L_{_fX}g)_{ij}=f(L_{_X}g)_{ij}+\left(\frac{\partial}{\partial x^i} f\right)g\left(X,\frac{\partial}{\partial x^j}\right)
    +\left(\frac{\partial}{\partial x^j} f\right)g\left(X,\frac{\partial}{\partial x^i}\right).$$

For our purposes, we may in addition assume that the curvature and the Ricci tensors of the initial metric admit a local $L^2$
and $L^p$ norm bounds respectively,  where $p>2$. By the argument in the next section, we can then show that the curvature and its
covariant derivative satisfy a local heat equation.  Also, they can be shown to satisfy $L^2$ energy bounds that are independent
of $\varepsilon >0$.  So (\ref{3.2}) has a solution for some time interval $[0,T)$, where $T$ is independent of $\varepsilon$.
Thus as $\varepsilon\rightarrow 0$, the solution of (\ref{3.2}) converges to a solution of (\ref{3.1}).

\vskip 1cm

\section{Smoothing a Riemannian Metric}
\ \ \ Let $X$ be a smooth manifold with Riemannian metric $g_0$ and $\Omega$ an open subset of $X$.
Let $\phi$ be a nonnegative smooth compactly supported function on $\Omega$.  Consider the
following evolution equation
\begin{align}\label{8}
    \left\{
    \begin{array}{ll}
       \dfrac{\partial g}{\partial t}&=-2\phi^2 \Ric(g),\\
       g(0)&=g_0.\\
    \end{array}
    \right.
\end{align}

It is easy to check that the curvature tensor $\Rm$ and Ricci tensor $\Ric$ satisfy the following equations respectively,
\begin{align*}
\frac{\partial \Rm}{\partial t}=\phi^2 (\triangle \Rm+Q(\Rm,\Rm))+2\phi a(\nabla \phi,\nabla \Rm)+b(\nabla \phi,\nabla \phi,\Rm)
        +\phi c(\nabla^2 \phi, \Rm)
\end{align*}
and
\begin{align*}
\frac{\partial \Ric}{\partial t}=\phi^2 (\triangle \Ric+Q(\Rm,\Ric))+2\phi a(\nabla \phi,\nabla \Ric)+b(\nabla \phi,\nabla \phi,\Ric)
        +\phi c(\nabla^2 \phi, \Ric).
\end{align*}
Notice that $\phi \in C_0^{\infty}(\Omega)$, we then have constant $c_1,\ c_2,\ ,c_3>0$ such that
$$\phi c(\nabla ^2 \phi,\Rm)\le -c_1 \phi \triangle \phi |\Rm|+c_2|\nabla\phi|^2|\Rm|+c_3\phi|\nabla \phi||\Rm|$$
and
$$\phi c(\nabla ^2 \phi,\Ric)\le -c_1 \phi \triangle \phi |\Ric|+c_2|\nabla\phi|^2|\Ric|+c_3\phi|\nabla \phi||\Ric|.$$
Then a direct calculation gives
\begin{align}\label{9}
\frac{\partial |\Rm|}{\partial t}\le \phi^2(\triangle |\Rm|+c_0|\Rm|^2)
    +2a\phi|\nabla \phi||\nabla \Rm|+b(|\nabla \phi|^2-\phi\triangle \phi)|\Rm|,
\end{align}
and
\begin{align}\label{10}
\frac{\partial |\Ric|}{\partial t}\le \phi^2(\triangle |\Ric|+c_0|\Rm||\Ric|)
    +2a\phi|\nabla \phi||\nabla \Ric|+b(|\nabla \phi|^2-\phi\triangle \phi)|\Ric|.
\end{align}

Again the results in this section are due to D. Yang \cite{4}.

\vskip 0.2cm
\begin{thm}\label{thm4}
There exist constant $C_1$ and $C_2$ such that if
$$\left(\int_{\Omega}|\Rm(g_0)|^2\der V_{g_0}\right)^{\frac 12}\le [\ C_1C_s(\Omega)\ ]^{-1}$$
and for any $p>2$,
$$\left(\int_{\Omega}|\Ric(g_0)|^p\der V_{g_0}\right)^{\frac 1p}<K,$$
then the equation (\ref{8}) has a smooth solution for $t\in[0,T)$,
where
$$T\ge \min \left(\|\nabla \phi\|^{-2}_{\infty}, C_2K^{-\frac{p}{p-2}}C_s(\Omega)^{-\frac{2}{p-2}} \right).$$
Moreover, for $t\in(0,T)$, the Riemannian curvature tensor satisfies the following bound,
\begin{align}\label{11}
\|\phi^2 \Rm\|_{\infty}\le C_3 C_s(\Omega)(t\|\nabla \phi\|^2_{\infty}+1)t^{-1}.
\end{align}
Here $C_1$ and $C_3$ only depend on the dimension of $X$; $C_2$ depends on the dimension of $X$ and $p$.
\end{thm}
\noindent\textit{Proof: } By Theorem \ref{thm3}, the equation (\ref{8}) has a smooth solution on a sufficiently
small time interval starting at $t=0$.
Let $[0,T_{\max})$ be a maximal time interval on which (\ref{8}) has a smooth solution and such that the
following hold for each metric $g(t)$,
\begin{align}
\label{12} \|f\|^2_{\psi}\le \psi A_0 \|\nabla f\|^2_2,&\ f\in C_0^{\infty}(\Omega);\\
\label{13} \frac 12 g_0 \le g(t) &\le 2 g_0;\\
\label{14} \|\Rm (g(t))\|_2\le& 2(C_1A_0)^{-1}.
\end{align}
Suppose that $T_{\max}<T_0=\min(\|\nabla \phi\|^{-2}_{\infty},\ C_2K^{-\frac{p}{p-2}}A^{-\frac {2}{p-2}})$.
We will show that this leads to a contradiction.

First, notice that the curvature tensor $\Rm$ satisfies (\ref{9}),
then according to the proof of Theorem \ref{thm2}, we have
\begin{align*}
\|\Rm(g(t))\|_2 < &e\|\Rm(g_0)\|_2\le 2e[C(n)4 A_0]^{-1}< 2 [C(n)A_0]^{-1},
\end{align*}
which implies a strict inequality for (\ref 14).

Next, since the Ricci curvature Ric satisfies (\ref{10}), then Corollary \ref{cor1} implies that
$$|\phi^2\Ric(g(t))|\le C_2A_0^{\frac 2p}(1+t\|\nabla \phi\|^2_{\infty})^2t^{-\frac{2}{p}}K.$$
Applying the bound on $\Ric$ to the following
$$\left|\frac{\der}{\der t}\int f^p\der V_g\right|\le 2\|\phi^2\Ric\|_{\infty}\int f^p\der V_g,$$
we have
$$\log\frac{\|f\|_p(t)}{\|f\|_p(0)}<\log 2.$$
The differential inequality
$$\left|\frac{\der}{\der t}\int |\nabla f|^2\der V_g\right|\le 2\|\Ric\|_{\infty}\int |\nabla f|^2\der V_g$$
leads to an analogous estimate.  Therefore, it follows that for any $t\le T_0$,
$$\|f\|^2_4(t)<2\|f\|^2_4(0)\le 2A_0\|\nabla f\|^2_2(0)<4A_0\|\nabla f\|^2_2(t),$$
that is to say (\ref{12}) holds with strict inequality.

To show that (\ref{13}) holds with strict inequality, we use Hamilton's trick.  Simply fix a tangent vector $v$
with respect to $g(t)$, then
$$\frac{\der}{\der t}|v|^2_{g(t)}=\frac{\der}{\der t}(g_{ij}(t)v^i v^j)=g_{ij}^{\prime}(t)v^i v^j$$
implies
$$\left|\frac{\der}{\der t}\log |v|^2_{g(t)}\right|\le |g_{ij}^{\prime}(t)|\le 2 \phi^2 |\Ric|.$$
So for $0\le t\le T_2 < T_0$,
$$\log \frac{|v|^2_{g(t)}}{|v|^2_{g(0)}}\le \int_0^{T_2}|g_{ij}^{\prime}(t)|\der t\le 2\|\phi^2\Ric\|_{\infty}T_2<\log 2,$$
which implies
$$\frac{1}{2}|v|^2_{g(0)}<|v|^2_{g(t)}<2|v|^2_{g(0)},$$
for $t<T_0$.

Finally, by differentiating the evolution equation for $\Rm$, we see that the covariant derivatives of $\Rm$ satisfy
evolution equations for which $L^2$ energy bounds can be obtained.
Therefore we can use Hamilton's argument in \S 14 of \cite{h} to show that $g(t)$ has a smooth limit as $t\rightarrow T_{\max}$.
If $T_{\max}<T_0$, we would be able to extend the solution to (\ref 8) smoothly beyond $T_{\max}$ with (\ref{12}), (\ref{13}) and (\ref{14})
still holding.  This contradicts the assumed maximality of $T_{\max}$. Hence, we conclude that $T_{\max}\ge T_0$.

The estimate (\ref{11}) follows from Theorem \ref{thm2}. $\Box$
\vskip 1cm

\section{Local Volume Estimate}

\ \ \ We consider more generally any system of the type
\begin{equation}\label{15}
\triangle \Ric=\Rm*\Ric+\B,
\end{equation}
where $\B$ is the Bach tensor.  Recall that
$$B_{ij}=2\nabla^k\nabla^l W^+_{ikjl}+R^{kl} W^+_{ikjl}.$$
We assume the following local Sobolev inequality,
$$\|f\|^2_{L^4(\Omega)}\le C_s(\Omega)\|\nabla f\|^2_2,\ \forall f\in C_0^{\infty}(\Omega).$$

\vskip 0.2cm
\begin{lem}\label{lem5}
There exist constant $\varepsilon,\ C$ such that if $\|\Rm\|_{L^2(B(p,r))}\le \varepsilon$ and
$\B\in L^2(B(p,r))$, then
$$\left\{\int_{B(p,\frac{r}{2})}|\Ric|^4\der V_g\right\}^{\frac 12}\le
    \frac{C}{r^2}\left(\int_{B(p,r)}|\Ric|^2\der V_g\right)+C\left(\int_{B(p,r)}|\Ric|^2\right)^{\frac 12}
    \left(\int_{B(p,r)}|B|^2\right)^{\frac 12}.$$
\end{lem}
\noindent \textit{Proof: }
From (\ref{15}), it follows that
$$\triangle |\Ric|\ge -|\Rm||\Ric|-|\B|.$$
We may assume that $r=1$.  The lemma then follows by scaling the metric.
Let $0\le \phi\le 1$ be a function supported in $B(p,1)$, then
\begin{align*}
\int_{B(p,1)}\phi^2|\Ric|^2|\Rm|\ge&\int \phi^2|\Ric|(-\triangle |\Ric|-\B)\\
    =&\int \nabla(\phi^2|\Ric|)\cdot \nabla |\Ric|-\int \phi^2|\Ric||\B|\\
    \ge& -\delta^{-1}\int|\nabla \phi|^2 |\Ric|^2+(1-\delta)\int|\phi\nabla|\Ric||^2\\
        &-\left(\int\phi^2|\Ric|^2\right)^{\frac 12}\left(\int \phi^2 |\B|^2\right)^{\frac 12}.
\end{align*}

Next, using the Sobolev constant bound, we have
\begin{align*}
\left(\int(\phi|\Ric|)^4\right)^{\frac 12}\le C\int |\nabla (\phi|\Ric|)|^2
   \le C\int |\nabla \phi|^2|\Ric|^2+C\int \phi^2|\nabla |\Ric||^2.
\end{align*}
Choosing $\delta$ sufficiently small yields
\begin{align*}
\left(\int(\phi|\Ric|)^4\right)^{\frac 12}\le& C\int \phi^2|\Ric|^2|\Rm|+C\int|\nabla \phi|^2|\Ric|^2\\
        &+C\left(\int \phi^2|\Ric|^2\right)^{\frac 12}\left(\int \phi^2|\B|^2\right)^{\frac 12}\\
    \le & C\left(\int \phi^2|\Rm|^2\right)^{\frac 12}\left(\int\phi^2|\Ric|^4\right)^{\frac 12}
        +C\int|\nabla \phi|^2|\Ric|^2\\
        &+C\left(\int \phi^2|\Ric|^2\right)^{\frac 12}\left(\int \phi^2|\B|^2\right)^{\frac 12}.
\end{align*}
Therefore, for $\varepsilon$ sufficiently small, we have
$$\left(\int\phi^2|\Ric|^4\right)^{\frac 12}\le C\int |\Ric|^2
        +C\left(\int |\Ric|^2\right)^{\frac 12}\left(\int |\B|^2\right)^{\frac 12}.$$
We then choose the cut-off function $\phi$ such that $\phi\equiv 1$ in $B(p,\frac{1}{2})$,
$\phi=0$ for $r=1$, $|\nabla \phi|\le C$, and we have
$$\left(\int_{B(p,\frac 12)}|\Ric|^4\right)^{\frac 12}\le C\int_{B(p,1)}|\Ric|^2
+C\left(\int |\Ric|^2\right)^{\frac 12}\left(\int |\B|^2\right)^{\frac 12}.$$
Scaling the metric, we obtain the lemma. $\Box$

\vskip 0.2cm
\begin{lem}\label{lem6}
With the same assumption of Lemma \ref{lem5}, we have
\begin{align*}
\left(\int_{B(p,\frac{r}{4})}|\Rm|^4\der V_g\right)^{\frac 12}\le&
    C \left(\int_{B(p,r)}|\Ric|^2\right)^{\frac 12}
    \left(\int_{B(p,r)}|\B|^2\right)^{\frac 12}\\
    &+\frac{C}{r^2}\int_{B(p,r)}|\Rm|^2.
\end{align*}
\end{lem}
\noindent \textit{Proof: }
Again we may assume that $r=1$.  Let $\phi$ be a cut-off function in $B(p,1)$, such that,
$\phi\equiv 1$ in $B(p,\frac 12)$ and $|\nabla \phi|\le C$.
We have, by lemma \ref{lem5},
\begin{align*}
\int_{B(p,1)}\phi^2|\nabla \Ric|^2=&-\int \phi^2\langle \triangle \Ric, \Ric \rangle-2\int \phi \langle \nabla \Ric,\nabla \phi \cdot\Ric \rangle\\
    =&-\int \phi^2\langle \Rm*\Ric, \Ric \rangle-\int \phi^2\langle \B, \Ric \rangle
        -2\int \phi \langle \nabla \Ric,\nabla \phi \cdot\Ric \rangle\\
    \le& C\left(\int \phi^2|\Rm|^2\right)^{\frac 12}\left\{\int_{B(p,1)}|\Ric|^2
        +\left(\int |\Ric|^2\right)^{\frac 12}\left(\int |\B|^2\right)^{\frac 12}\right\}\\
    &+C\int |\Ric|^2+C\delta\int \phi^2|\nabla\Ric|^2+C\left(\int |\B|^2\right)^{\frac 12}\left(\int |\Ric|^2\right)^{\frac 12}.
\end{align*}
By choosing $\delta$ small and $\varepsilon<1$, we have
\begin{align}\label{16}
\int_{B(p,\frac 12)}|\nabla\Ric|^2\le&(1+\varepsilon)C\int_{B(p,1)}|\Ric|^2 \nonumber\\
        &+(1+\varepsilon)C\left(\int_{B(p,1)}|\Ric|^2\right)^{\frac 12}\left(\int_{B(p,1)}|\B|^2\right)^{\frac 12}\nonumber\\
    \le & 2C\int_{B(p,1)}|\Ric|^2\nonumber\\
        &+2C\left(\int_{B(p,1)}|\Ric|^2\right)^{\frac 12}\left(\int_{B(p,1)}|\B|^2\right)^{\frac 12}.
\end{align}
Next, let $\phi$ be a cutoff function in $B(p, \frac 12)$, such that $\phi \equiv 1$ in $B(p,\frac 14)$ and $|\nabla \phi|\le C$.
Recall that
$$\triangle \Rm=L(\nabla^2 \Ric)+\Rm*\Rm,$$
where $L(\nabla^2\Ric)$ denotes a linear expression in second derivatives of the Ricci tensor.
We then have
\begin{align*}
\int_{B(p,\frac 12)}\langle\triangle\Rm,\phi^2\Rm\rangle=& \int\langle\nabla^2\Ric+\Rm*\Rm,\phi^2\Rm\rangle\\
    =&-\int\langle 2\phi\nabla \Ric,(\nabla \phi)\Rm\rangle-\int \phi^2\langle\nabla\Ric,\nabla\Rm\rangle\\
        &+\int \langle\Rm*\Rm,\phi^2\Rm\rangle.
\end{align*}
This yields
\begin{align*}
\left|\int_{B(p,\frac 12)}\langle\triangle\Rm,\phi^2\Rm\rangle\right|\le&
        C\int\phi^2|\nabla \Ric|^2+C\int|\nabla \phi|^2\nabla \Rm|^2\\
        &+\frac{C}{\delta}\int\phi^2|\nabla \Ric|^2+C\delta\int\phi^2|\nabla \Rm|^2+C\int\phi^2|\Rm|^3.
\end{align*}
Integrating by parts,
\begin{align*}
\int_{B(p,\frac 12)}\phi^2|\nabla \Rm|^2=&\int\langle 2 \phi\nabla\Rm,(\nabla \phi)\Rm\rangle
        -\int\phi^2\langle \triangle\Rm,\Rm \rangle\\
        \le& \frac{C}{\delta}\int|\nabla\phi|^2|\Rm|^2+2C\delta \int\phi^2|\nabla\Rm|^2  +C\int|\nabla\phi|^2|\Rm|^2\\
        &+\frac{C'}{\delta}\int\phi^2|\nabla\Ric|^2+C\int\phi^2|\Rm|^3.
\end{align*}
Choosing $\delta$ sufficiently small and using (\ref{16}), we obtain
\begin{align*}
\int\phi^2|\nabla \Rm|^2
    \le& C\int |\nabla\phi|^2|\Rm|^2 +C\int|\Rm|^2+C\int\phi^2|\Rm|^3\\
    &+C\left(\int |\Ric|^2\right)^{\frac 12}\left(\int |\B|^2\right)^{\frac 12}.
\end{align*}
Using the Sobolev inequality,
\begin{align*}
\left(\int|\phi\Rm|^4\right)^{\frac 12}\le&C\int|\nabla|\phi\Rm||^2\\
    \le&C\int|\nabla \phi|^2|\Rm|^2+C\int \phi^2|\nabla|\Rm||^2\\
    \le&C\int|\Rm|^2+C\left(\int|\Rm|^2\right)^{\frac 12}\left(\int\phi^4|\Rm|^4\right)^{\frac 12}\\
        &+C\left(\int |\Ric|^2\right)^{\frac 12}\left(\int |\B|^2\right)^{\frac 12}.
\end{align*}
Therefore by choosing $\varepsilon$ small, we obtain
\begin{align*}
\left(\int_{B(p,\frac 14)}|\Rm|^4\right)^{\frac 12}
    \le C\int|\Rm|^2+C\left(\int |\Ric|^2\right)^{\frac 12}\left(\int |\B|^2\right)^{\frac 12}.
\end{align*}
Scaling the metric, we obtain the lemma.  $\Box$

\vskip 0.2cm
\begin{thm}\label{thm5}
Assume that (\ref{15}) is satisfied.  Let $B(p,r)$ be a geodesic ball around the point $p$.
Then there exist constants $\varepsilon_0$, $C$ (depending on the Sobolev constant $C_s(B(p,r))$) such that if
$$\|\Rm\|_{L^2(B(p,2r))}\le \varepsilon_0$$
and $\B\in L^2(B(p,2r))$, then
$$\Vol(B(p,r))\le C r^4.$$
\end{thm}
\noindent \textit{Proof: } We assume that $r=1$.  By lemma
\ref{lem5} and \ref{lem6}, we have
$$\int |\Rm|^3\le \left(\int|\Rm|^2\right)^{\frac{1}{2}}\left(\int|\Rm|^4\right)^{\frac{1}{2}}\le C.$$
Then for $\varepsilon_0$ suitably chosen, by Theorem \ref{thm4}, the local Ricci flow
\begin{align*}
    \left\{
        \begin{array}{ll}
        \dfrac{\partial g(t)}{\partial t}&=-2\phi^2\Ric(g(t)),\\
        g(0)&=g
        \end{array}
    \right.
\end{align*}
has a smooth solution for $t\in [0,T)$, where
$$T\ge \min(\|\nabla \phi\|^{-2}_{\infty},C^{-3}C_s^{-2})$$
and for $t\in(0,T)$, the Riemannian curvature tensor satisfies the following bound
$$\|\phi^2\Rm\|_{\infty}\le CC_s(t\|\nabla \phi\|^2_{\infty}+1)t^{-1}.$$
Therefore,
$$\Vol_{g(t)}(B(p,1))\le C$$
for any fixed $t\in (0,T)$.  Since we can find a constant $C$ such that
$$\frac{1}{C}g\le g(t)\le Cg,$$
then for the metric $g$ we still obtain the volume estimate
$$\Vol_g(B(p,1))\le C.$$
This proves the theorem.  $\Box$

\end{document}